\documentstyle[12pt]{article}
\begin{document}
\begin{center}
{\Large \bf
Vertex Operator of ${U_q}(\widehat{B_l})$ for Level One
\par}
\lineskip .75em
\vskip 3em

\large{Bai-Qi Jin \footnote{e-mail address::jinbq@bepc3.ihep.ac.cn}}  \\
{\footnotesize Institute of High Energy Physics, P.O.Box 918(4), Beijing
100039, P. R. China}\\
\vspace{1.5em}

\large{Shan-You Zhou}  \\
{\footnotesize Department of Mathematics, Graduate School, University of
Science and Technology of China\\
P.O.Box 3908,Beijing 100039,P.R.China}

\vskip 3.0em

{\bf Abstract}
\end{center}
\par
In this paper,we give the explicit formulae of vertex operators of
${U_q}(\widehat{B}_l)$ for level-one as operators on the Fock space.
Meanwhile, we point out that the free field realization (by one
fermionic field and $l$ bosonic fields)
of highest weight module with highest weight $\Lambda_l$ has two
irreducible modules.

\vspace{1.0cm}
\noindent
{\bf
1 Introduction}

One of the central subjects of mathematical physics has been the studies
on exactly solvable models in two dimensions for many years.
The central problem is eigenstate and correlation functions (form-factors)
in exactly solvable models. In \cite{DFJM} \cite{JM} a new scheme was
given  for solving the six-vertex model and associated XXZ chain in the
antiferromagnetic regime  using the newly discovered quantum affine symmetry
of the system. The approach of that paper has been extended to higher
spin-chains \cite{Konno} \cite{BW} \cite{IIJM} \cite{Idzumi},to the
higher rank case \cite{YK} \cite{DO}
and to the ABF models \cite{JMO}. The analogous approach in integrable
massive field theory also has been developed by Lukyanov \cite{Luky}.
All of these papers are concerned with models constructed on the quantum
affine algebra $U_q(\widehat{A_l})$.
The key object in this approach is
vertex operators which have first been introduced by Frenkel and Reshetikhin
\cite{FR}. In order to extend that scheme to the models with symmetry
of $U_q(\widehat{B_l})$,
we construct the $q$-vertex operators related to $U_q(\widehat{B_l})$
in this paper.
The physical application will be developed else where.In section $2$,
we briefly recall the free field realization of $U_q(\widehat{B_l})$.
In section $3$, we deduce the vertex operators.

\vspace{1.0em}
\noindent
{\bf
2 Free Field Realization of $U_q(\widehat{B}_l)$}

In this section we briefly recall the free field realization of
$U_q(\widehat{B_l})$ \cite{Bernard}

\noindent
$2.1$ Definition of $U_q(\widehat{B_l})$.

First, fix some notations we will use. Let $\widehat{P}={\bf Z} {\Lambda}_0
\oplus {\bf Z} {\Lambda}_1 \oplus \cdots \oplus {\bf Z} {\Lambda}_l
\oplus {\bf Z} {\delta}$ be the weight lattice of $U_q(\widehat{B_l})$
and let $\alpha_i=\displaystyle \sum_{j=0}^l a_{ji} {\Lambda}_j +
{\delta}_{i,0} \delta$
($i=0,1,\cdots,l$) be simple roots.So we have $\delta=\alpha_0+
\alpha_1+2(\alpha_2+\cdots +\alpha_l)$ and $\widehat{Q}={\bf Z} \alpha_0
\oplus {\bf Z} \alpha_1 \oplus \cdots \oplus {\bf Z} \alpha_l$ being its
root lattice.The symmetric bilinear form (,) on $\widehat{P}$ is defined by
$$({\Lambda}_0,{\Lambda}_i)=0~~~i=0,1,\cdots,l ,~~~({\Lambda}_0,\delta)=1$$
$$({\Lambda}_i,{\alpha}_j)=d_i {\delta}_{ij}$$
we can get
$$({\Lambda}_1,\delta)=1,~~~(\delta,\delta)=0,~~~({\Lambda}_l,\delta)=1$$
$$({\Lambda}_i,{\Lambda}_j)=
d_i ({\overline{A}}^{-1})_{ij}~~~1\leq i,j\leq l,~~~({\Lambda}_i,\delta)=2
{}~~~1<i<l$$
$$({\Lambda}_i,\delta)=2~~1\leq i\leq l,~~~(\alpha_i,\alpha_j)=d_i a_{ij}$$
where $A=(a_{ij})_{i,j=0}^l$ is the Cartan matrix of $U_q(\widehat{B_l})$ and
$\overline{A}$ is the matrix of $A$ removed by first line and first
column.

For $U_q(\widehat B_2)$,$$d_0=d_1=2 d_2=1,$$
$$A=\left[
\begin{array}{ccc}
2&0&-1\\
0&2&-1\\
-2&-2&2
\end{array}
\right]$$
and for $U_q(\widehat B_l)$ $(l\geq 3)$,$$d_0=d_1= \cdots=d_{l-1}=2 d_l=1,$$
$$A=\left[
\begin{array}{ccccccc}
2&0&-1&0&\cdots&0&0\\
0&2&-1&0&\cdots&0&0\\
-1&-1&2&-1&\cdots&0&0\\
0&0&-1&2&\cdots&0&0\\
\cdot&\cdot&\cdot&\cdot&\cdots&\cdot&\cdot\\
0&0&0&0&\cdots&2&-1\\
0&0&0&0&\cdots&-2&2
\end{array}
\right]$$
Define the dual space of $\widehat{P}$ as $\widehat{P^*}={\bf Z} h_0
\oplus {\bf Z} h_1 \oplus \cdots \oplus {\bf Z} h_l \oplus {\bf Z} d$.
The dual pairing $\langle,\rangle$ is defined by \\
$$\langle h_i, \lambda \rangle: =d_i^{-1} ({\alpha}_i,\lambda)~~~\lambda
\in \widehat{P}$$
$$\langle d,\lambda \rangle :=({\Lambda}_0,\lambda)~~~\lambda \in \widehat{P}$$
Later we will use the weight lattice and the root lattice of $U_q(B_l)$.
The weight lattice $P={\bf Z} \lambda_1 \oplus {\bf Z} \lambda_2
\oplus \cdots \oplus {\bf Z} \lambda_l$,
$\lambda_i=\Lambda_i-a_i^{\vee}\Lambda_0$.
where $a_0^{\vee}=a_1^{\vee}=a_l^{\vee}=1$,$a_2^{\vee}=\cdots=a_{l-1}^{\vee}=2$
The root lattice of $U_q(B_l)$ is $Q={\bf Z} \alpha_1 \oplus
{\bf Z} {\alpha_2} \oplus \cdots \oplus
{\bf Z} {\alpha_i}$ .
In this article ,we assume $-1<q<0$.
and use the following standard notations:

$$\omega=\frac{1}{(q^{\frac{1}{2}}+q^{-\frac{1}{2}})},~~~q_i=q^{d_i}$$
$$[m]_{q_i}=\frac{{q_i}^m-{q_i}^{-m}}{q_i-{q_i}^{-1}}~~~~~m \in {\bf Z}$$.
$$\left[
\begin{array}{c}
n\\
m
\end{array}
\right]_{q_i}=\frac{[1]_{q_i} [2]_{q_i} \cdots [n]_{q_i}}{[1]_{q_i} [2]_{q_i}
\cdots [m]_{q_i} [1]_{q_i} [2]_{q_i} \cdots [n-m]_{q_i}}$$
and when $q_i=q$,we omit the index.
Quantum affine algebra $U_q(\widehat{B}_l)$ is an associative
algebra over ${\cal C}$ with unity
generated by $e_i$,$f_i$ ,$q_i^{\pm h_i}$ $(i=0,1,2,\cdots,l)$ and $q^{\pm d}$.
$U^{\prime}_q(\widehat{B}_l)$ is the subalgebra of
$U_q(\widehat{B}_l)$ generated by $\{ q_i^{\pm h_i},e_i,f_i
\mid 0\leq i\leq l\}$.
The defining relations are as follows ($\forall h,h^{\prime} \in
\widehat{P^*}$):
$$q^h q^{h'}=q^{h+h'},~~~q^0=1$$
$$q^h e_i q^{-h}=q^{\langle h,{\alpha}_i \rangle} e_i$$
$$q^h f_i q^{-h}=q^{- \langle h,{\alpha}_i \rangle} f_i$$
$$[e_i,f_j]={\delta}_{i,j} \frac{q_i^{h_i}-q_i^{-{h_i}}}{q_i-q_i^{-1}}$$
$$\sum_{n=0}^{1-a_{ij}} (-1)^n {\left[
\begin{array}{c}
1-a_{ij}\\
n
\end{array}
\right]}_{q_i} {e_i}^{1-a_{ij}-n} e_j {e_i}^n=0,~~~i\not= j$$
$$\sum_{n=0}^{1-a_{ij}} (-1)^n {\left[
\begin{array}{c}
1-a_{ij}\\
n
\end{array}
\right]}_{q_i} {f_i}^{1-a_{ij}-n} f_j {f_i}^n=0,~~~i\not= j$$
The algebra $U_q(\widehat{B_l})$ has a Hopf algebra structure
with the following coproduct $\Delta$ : $U_q(\widehat{B_l})\rightarrow
U_q(\widehat{B_l})  \otimes U_q(\widehat{B_l})$
$$\Delta (e_i)=e_i \otimes 1 +q_i^{h_i} \otimes e_i,~~~\Delta (f_i)=f_i
\otimes q_i^{-h_i} + 1\otimes f_i$$
$$\Delta (q^h) =q^h \otimes q^h,~~~ h\in P^*$$
As above this algebra is defined through Chevalley generators.

\noindent
$2.2$ Drinfeld realization \cite{Drinfeld}

For free field realization, $U_q(\widehat{B_l})$ defined by Chevalley
generators is not convenient.Drinfeld has given
another realization in 1988 \cite{Drinfeld}.The new realization is
an associative algebra generated by the elements $\{x_i^{\pm}(n),a_i(m),
\gamma^{\pm{\frac{1}{2}}},q_i^{h_i} \mid 1\leq i \leq l,m \in
{\bf Z}\backslash\{0\},n\in {\bf Z}\}$
satifying the follow relations:
\begin{center}
${\gamma}^{\pm\frac{1}{2}}$ is the center of the algebra.
\end{center}
$$[a_i(n),a_j(m)]={\delta}_{n+m,0} \frac{1}{n} [n a_{ij}]_{q_i}
\frac{{\gamma}^n - {\gamma}^{-n}}{q_j-{q_j^{-1}}}$$
$$q_i^{h_i} {x^{\pm}}_j(n) q_i^{-h_i} ={q_i}^{\pm a_{ij}} {x_j}^{\pm}(n)$$
$$[a_i(n),{x_j}^{\pm}(m)] = \pm \frac{1}{n} [n a_{ij}]_{q_i}
{\gamma}^{\mp \frac{\mid n \mid}{2}} {x_j}^{\pm}(n+m)$$
$${x_i}^{\pm}(n+1) {x_j}^{\pm}(m)-{q_i}^{\pm a_{ij}} {x_j}^{\pm}(m)
{x_i}^{\pm}(n+1)={q_i}^{\pm a_{ij}} {x_i}^{\pm}(n)  {x_j}^{\pm}(m+1)
-{x_j}^{\pm}(m+1){x_i}^{\pm}(n)$$
$$[{x_i}^+(n),{x_j}^{-}(m)]={\delta}_{i,j} \frac{{\gamma}^{\frac{n-m}{2}}
{{\varphi}_i}^+(n+m)-{\gamma}^{-\frac{n-m}{2}}
{{\varphi}_i}^{-}(n+m)}{q_i-{q_i}^{-1}}$$
$$\sum_{\pi \in \sum_p} \sum_{k=0}^p (-1)^k {\left[
\begin{array}{c}
p\\
k
\end{array}
\right]}_{q_i} {x_i}^{\pm}(r_{\pi (1)})
\cdots {x_i}^{\pm}(r_{\pi (k)}) {x_j}^{\pm}(s) {x_i}^{\pm}(r_{\pi (k+1)})
\cdots {x_j}^{\pm}(r_{\pi (p)})=0$$
if $i \not= j$, for sequences of integers $r_1$,$\cdots$, $r_p$.
where $p=1-a_{ij}$,
$\sum_p$ is the symmetric group on $p$ letters, and the
${{\varphi}_i}^{\pm}(r)$
are determined by equating powers of $u$ in the formal power series
$$\sum_{r=0}^{\infty} {{\varphi}_i}^{\pm}(r) u^{\mp r}=q_i^{\pm h_i}
exp(\pm (q_i-{q_i}^{-1})\sum_{n=1}^{\infty} a_i(\pm n) u^{\mp n})$$

Define maps ${w_i}^{\pm}$: $U_q(\widehat{g})\rightarrow U_q(\widehat{g})$ by
$${w_i}^{\pm}a={x_i}^{\pm}(0) a -q_i^{\pm h_i} a q_i^{\mp h_i} {x_i}^{\pm}(0)$$
then, the isomorphism between the two realization is
$$q^{h_0}={\gamma} q^{-h_1-2 h_2 -\cdots -2 h_{l-1} - h_l},
{}~~~e_i={x_i}^+(0),~~~f_i={x_i}^-(0),~~~for~~ i =1,2,\cdots,l$$
$$e_0={w^-}_2 {w^-}_3 \cdots {w^-}_l {w^-}_l \cdots {w^-}_2 {X_1}^-(1)
q^{h_0} {\gamma}^{-1}=e^{\prime}_0 q^{h_0} {\gamma}^{-1}$$
$$f_0=q^{2l-3} {\omega}^2 q^{-h_0} {\gamma} {w^+}_2 {w^+}_3 \cdots {w^+}_l
{w^+}_l \cdots {w^+}_2 {X_1}^+(-1)= q^{-h_0} {\gamma} f^{\prime}_0 $$

\noindent
$2.3$ Central Extension of Weight Lattice

It is well known that it is enough to consider only the central extension of
root lattice $Q$ to construct the representation of $U_q(\widehat{B_l})$.
But, for the construction of the vertex operators it is necessary to define
the central extention of weight lattice  $P$.

Define the group algebra $C[\overline{P}]$ generated by the symbols
$\{e^{{\lambda}_1}$,$e^{{{\alpha}}_1}$,
$e^{{{\alpha}}_2}$,$\cdots$,$e^{{{\alpha}}_{l-1}}\}$  and satisfying
the following relations:
$$e^{{{\alpha}}_i} e^{{{\alpha}}_j}= (-1)^{({{\alpha}}_i,{{\alpha}}_j)}
e^{{{\alpha}}_j} e^{{{\alpha}}_i}
{}~~~~ 1\leq i,j\leq l-1$$
$$e^{{{\alpha}}_i} e^{{\lambda}_1}=(-1)^{\delta_{i,1}} e^{{\lambda}_1}
e^{{{\alpha}}_i} ~~~~ 1 \leq i\leq l-1$$
For $\alpha=m_0 {\lambda}_1 + m_1 {{\alpha}}_1 +m_2 {{\alpha}}_2+
\cdots +m_{l-1} {{\alpha}}_{l-1}$,we denote
$e^{\alpha}=e^{m_0 {\lambda}_1} e^{m_1 {{\alpha}}_1} e^{m_2 {{\alpha}}_2}
\cdots e^{m_{l-1} {{\alpha}}_{l-1}}$
So a simple calculation shows $$e^{{{\alpha}}_l} e^{{\lambda}_1}=-
e^{{\lambda}_1} e^{{{\alpha}}_l}$$
$$e^{{{\alpha}}_i} e^{{{\alpha}}_j}= (-1)^{({{\alpha}}_i,{{\alpha}}_j)}
e^{{{\alpha}}_j} e^{{{\alpha}}_i} ~~~~ 1\leq i,j\leq l$$

\noindent
$2.4$ Level One Module

Bernard \cite{Bernard} has given the free field realization
of level one module by using one fermionic field and $l$ bosonic fields.
Here,we reconstruct the three level one module and point out a fact having
been overlooked.

\noindent
Let
$$ H_i:= C[a_j(-m), {\Psi}(-k)
(1 \leq j \leq l,m\in {\bf Z}_{\geq 0},k \in {{\bf Z}_{\geq 0}+
\frac{1}{2}})]\otimes
C[\overline{Q}] e^{{\Lambda}_i}$$
for $i=0,1$, in ${\cal NS}$ cases.
$$ H_l:= C[a_j(-m), {\Psi}(-k) (1 \leq j \leq l,m\in {{\bf Z}_{\geq 0}},
k \in {\bf Z}_{\geq 0})]\otimes C[\overline{Q}] e^{{\Lambda}_i}$$
in ${\cal R}$ cases.

We define the operators $a_i(k)$ $(1 \leq i \leq l)$, ${\Psi}(k)$,
${\partial}_{{{\alpha}}_i}$, $e^{{{\alpha}}_i}$,
$d$ on $H_i$ $(i=0,1,l)$ as follows:

for $f\otimes e^{\beta}=a_{i_1}(-n_1) \cdots a_{i_k}(-n_k)
{\Psi}(-k_1)\cdots {\Psi}(-k_n) \otimes e^{\beta} \in H_i$,
$$ {a_j}(k). f \otimes e^{\beta}=\lbrace {
\begin{array}{cc}
{{a_j}(k)f \otimes e^{\beta}} & {(k<0)}\\
{[{a_j}(k),f] \otimes e^{\beta}} & {(k>0)}
\end{array}}$$
$$ {\Psi}(k). f \otimes e^{\beta}=\lbrace {
\begin{array}{cc}
{{\Psi}(k)f \otimes e^{\beta}}&{(k \leq 0)}\\
{\{{\Psi}(k),f\} \otimes e^{\beta}}&{(k>0)}
\end{array}}$$
$${\partial}_{{{\alpha}}_i}.f \otimes e^{\beta}= ({{\alpha}}_i,\beta) f
\otimes e^{\beta}$$
$$e^{{\alpha}_j}.f\otimes e^{\beta}=f \otimes e^{{\alpha}_j} e^{\beta}$$
$$d.f\otimes e^{\beta}= (-\sum_{l=1}^k n_l-\sum_{l=1}^n (k_l-\frac{1}{2})-
\frac{(\beta,\beta)}{2} +\frac{({\Lambda}_i,{\Lambda}_i )}{2}) f\otimes
e^{\beta}$$

\noindent
Let $$x^{\pm}_i(z):=\sum_{n \in {\bf Z}}x^{\pm}_i(n) z^{-n}~~~(1 \leq i
\leq l-1)$$

We define the action of $U_q(\widehat{B_l})$:
$$ \gamma \rightarrow q,~~~ q_j^{h_j} \rightarrow q^{{\partial}_{{\alpha_j}}},
{}~~~ (1 \leq j \leq l)$$
$${x_i}^{\pm}(z)\rightarrow z exp({\pm}\sum_{k=1}^{\infty}
{\frac{a_i(-k)}{[k]}} q^{\mp \frac{k}{2}} z^k) exp(\mp
\sum_{k=1}^{\infty} \frac{a_i(k)}{[k]} q^{\mp \frac{k}{2}} z^{-k})
e^{\pm {{\alpha}}_i} z^{\pm {\partial}_{{\alpha}_i}}$$
$${x_l}^{\pm}(z)\rightarrow z^{\frac{1}{2}} exp({\pm}\sum_{k=1}^{\infty}
{\frac{a_l(-k)}{[k]}} q^{\mp \frac{k}{2}} {\omega} z^k) exp(\mp
\sum_{k=1}^{\infty} \frac{a_l(k)}{[k]} q^{\mp \frac{k}{2}} {\omega} z^{-k})
e^{\pm {{\alpha}}_l} z^{\pm {\partial}_{{\alpha}_l}} {\Psi}(z)$$
where
$${\Psi}(z)=\sum_n {\Psi}(n) z^{-n}$$
with $n\in {\bf Z}+{\frac{1}{2}}$(or $n\in {\bf Z}$) in the ${\cal NS}$
(or ${\cal R}$) cases,respectively,
and $$\{{\Psi}(n),{\Psi}(m)\}=(q^n+q^{-n}){\delta}_{n+m,0}$$
$$[a_i(n),a_j(m)]={\delta}_{n+m,0} \frac{1}{n} [n a_{ij}]_{q_i}
\frac{q^n - q^{-n}}{q_j-{q_j^{-1}}}$$
$$[a_i(n),{\Psi}(m)]=0$$

\noindent
Define
$$G=(-1)^{2 \displaystyle\sum_{i=1}^l \partial_{\lambda_i}+N_F}$$
for ${\cal NS}$ cases.
$$G=(-1)^{2 \displaystyle \sum_{i=1}^l \partial_{\lambda_i}+N_F-2
\sum_{i=1}^l (\lambda_i,\lambda_l)}$$
for ${\cal R}$ cases.
where $N_F$ denotes the fermion's number operator.
We easily find $G$ commutes with all the elements of $U_q(\widehat{B_l})$.
Through the eigenvales of $G$,we can divide Fock space into irreducible ones.
Four irreducible $U_q(\widehat{B_l})$-modules $V(\Lambda_0)$,$V(\Lambda_1)$,
$V(\Lambda_l)$,$V({\Lambda_l}^{\prime})$ whose highest weight vectors are
$1\otimes e^{\Lambda_0}$,$1\otimes e^{\Lambda_1}$,$1\otimes e^{\Lambda_l}$,
$\Psi(0)\otimes e^{\Lambda_l}$,respectively.The first two modules are in
${\cal NS}$ cases and the others are in ${\cal R}$ cases.We have not
seen that the
reducibility in ${\cal R}$ cases has been discussed  before.

\vspace{1.0em}
\noindent
{\bf
3. Vertex Operators}

Let us review the definition and some properties of the vertex operators.

\noindent
$3.1$ Finite dimensional $U^{\prime}_q(\widehat{B_l})$ module

Let $V$ be a finite dimensional $U^{\prime}_q(\widehat{B_l})$-module
with basis $\{v_m \mid 1 \leq m \leq 2l+1\}$,
and the representation of $e_i,f_i,h_i$ is as follows:
$$e_i={f_i}^t=E_{i,i+1}+E_{2l-i+1,2l-i+2},~~~i\not= 0,l$$
$$e_l={f_l}^t={\omega}^{-\frac{1}{2}} (E_{l,l+1}+E_{l+1,l+2})$$
$$e_0={f_0}^t=E_{2l,1}+E_{2l+1,2}$$
$$h_i=E_{i,i}-E_{i+1,i+1}+E_{2l-i+1,2l-i+1}-E_{2l-i+2,2l-i+2},~~~i\not= 0,l$$
$$h_0=-E_{1,1}-E_{2,2}+E_{2l,2l}+E_{2l+1,2l+1}$$
$$h_l=2E_{l,l}-2E_{l+2,l+2}$$
Define the $U_q$-module structure on $V_z$ as follows.
$$e_i(v_m \otimes z^n)=e_i v_m \otimes z^{n+{\delta}_{i,0}},~~~f_i(v_m
\otimes z^n)=f_i v_m \otimes z^{n-{\delta}_{i,0}}$$
$$h_i(v_m \otimes z^n)=h_i v_m \otimes z^n,~~~d (v_m \otimes z^n)=n v_m
\otimes z^n$$
We call $V_z$ the affinization of $V$ as a $U_q(\widehat{B_l})$-module
of level zero.

\noindent
$3.2$ Definition of $q$-Vertex Operators

The intertwiners of $U_q(\widehat{B_l})$-modules
$${{\widehat{\Psi}}_{\Lambda_i}}^{{\Lambda_j}V}(z) :
V(\Lambda_i) \rightarrow V(\Lambda_j) \otimes V_z$$
are called type-{\rm I} vertex operators.and the operators
$${{\widehat{\Psi}}_{\Lambda_i}}^{V{\Lambda_j}}(z) :V(\Lambda_i)
\rightarrow V_z\otimes V(\Lambda_j) $$
are called type-{\rm II} vertex operators.where $\otimes$ is the tensor product
 with an appropriate completion.
Denote the vertex operators as a formal series
$${{\widehat{\Psi}}_{\Lambda_i}}^{{\Lambda_j}V}(z)=
\sum_{m=1}^{2l+1} {{\widehat{\Psi}}_{{\Lambda_i}m}}^{{\Lambda_j}V}(z)
\otimes v_m$$
$${{\widehat{\Psi}}_{\Lambda_i}}^{V{\Lambda_j}}(z)=
\sum_{m=1}^{2l+1} v_m\otimes {{\widehat{\Psi}}_{{\Lambda_i}m}}^{V
{\Lambda_j}}(z)$$
$${{\widehat{\Psi}}_{{\Lambda_i}m}}^{{\Lambda_j}V}(z)=\sum_{n\in
{\bf Z}}{{\widehat{\Psi}}_{{\Lambda_i}m}}^{{\Lambda_j}V}(n) z^n$$
$${{\widehat{\Psi}}_{{\Lambda_i}m}}^{V{\Lambda_j}}(z)=\sum_{n\in
{\bf Z}}{{\widehat{\Psi}}_{{\Lambda_i}m}}^{V{\Lambda_j}}(n) z^n$$
There exist four type-{\rm I} (respectively type-{\rm II})
vertex operators \cite{FR} \cite{DJO}.
$${{\widehat{\Psi}}_{\Lambda_0}}^{{\Lambda_1}V}(z) :
V(\Lambda_0) \rightarrow V(\Lambda_1) \otimes V_z$$
$${{\widehat{\Psi}}_{\Lambda_1}}^{{\Lambda_0}V}(z) :
V(\Lambda_1) \rightarrow V(\Lambda_0) \otimes V_z$$
$${{\widehat{\Psi}}_{\Lambda_l}}^{{\Lambda_l^{\prime}}V}(z) :
V(\Lambda_l) \rightarrow V(\Lambda_l^{\prime}) \otimes V_z$$
$${{\widehat{\Psi}}_{\Lambda_l^{\prime}}}^{{\Lambda_l}V}(z) :
V(\Lambda_l^{\prime}) \rightarrow V(\Lambda_l) \otimes V_z$$
We can impose the nomalization condition.\\
{}~~~${{\widehat{\Psi}}_{{\Lambda_0}}}^{{\Lambda_1}V}(z).(1
\otimes e^{\Lambda_0}) =(1\otimes e^{\Lambda_1})
\otimes v_{2l+1}+$(terms of positive powers in $z$)\\
{}~~~${{\widehat{\Psi}}_{{\Lambda_1}}}^{{\Lambda_0}V}(z).(1
\otimes e^{\Lambda_1}) =(1\otimes e^{\Lambda_0})
\otimes v_1+$(terms of positive powers in $z$)\\
{}~~~${{\widehat{\Psi}}_{{\Lambda_l}}}^{{\Lambda_l^{\prime}}V}(z).(
1\otimes e^{\Lambda_l})=(\Psi (0)\otimes e^{\Lambda_l})\otimes
v_{l+1}+$(terms of positive powers in $z$)\\
{}~~~${{\widehat{\Psi}}_{{\Lambda_l^{\prime}}}}^{{\Lambda_l}V}(z).(\Psi (0)
\otimes e^{\Lambda_l})=(1\otimes e^{\Lambda_l})\otimes
v_{l+1}+$(terms of positive powers in $z$)\\
Later we will find the formulae of
${{\widehat{\Psi}}_{{\Lambda_l}}}^{{\Lambda_l^{\prime}}V}(z)$
and ${{\widehat{\Psi}}_{{\Lambda_l^{\prime}}}}^{{\Lambda_l}V}(z)$
identical.So we will
denote them as ${{\widehat{\Psi}}_{{\Lambda_l}}}^{{\Lambda_l}V}(z)$.

\noindent
$3.3$ Vertex Operator of Type-{\rm I}

{}From the intertwining relation (we only need to describe the type-{\rm I}
vertex operator.the type-{\rm II} is quite parallel).
$${\Delta}(x) \circ {{\widehat{\Psi}}_{\Lambda_i}}^{{\Lambda_j}V}(z)=
{{\widehat{\Psi}}_{\Lambda_i}}^{{\Lambda_j}V}(z) \circ x $$
then we can get the commutation relations between the type-{\rm I} vertex
operator with Chevalley generators.
Here we write the relations partially.
$$q f_2 {{\widehat{\Psi}}_{{\Lambda_{\sigma_1}}3}}^{{\Lambda_{\sigma_2}}V}(z)
+ {{\widehat{\Psi}}_{{\Lambda_{\sigma_1}}2}}^{{\Lambda_{\sigma_2}}V} (z)=
  {{\widehat{\Psi}}_{{\Lambda_{\sigma_1}}3}}^{{\Lambda_{\sigma_2}}V}(z) f_2$$
  $$q f_3 {{\widehat{\Psi}}_{{\Lambda_{\sigma_1}}4}}^{{\Lambda_{\sigma_2}}V}(z)
  +{{\widehat{\Psi}}_{{\Lambda_{\sigma_1}}3}}^{{\Lambda_{\sigma_2}}V} (z)
  ={{\widehat{\Psi}}_{{\Lambda_{\sigma_1}}4}}^{{\Lambda_{\sigma_2}}V} (z) f_3$$
  $$ \vdots ~~~~\vdots~~~~\vdots$$
$$q f_{l-1}
{{\widehat{\Psi}}_{{\Lambda_{\sigma_1}}l}}^{{\Lambda_{\sigma_2}}V}(z)+
  {{\widehat{\Psi}}_{{\Lambda_{\sigma_1}}{l-1}}}^{{\Lambda_{\sigma_2}}V}(z)
  ={{\widehat{\Psi}}_{{\Lambda_{\sigma_1}}l}}^{{\Lambda_{\sigma_2}}V} (z)
f_{l-1}$$
  $$f_l
{{\widehat{\Psi}}_{{\Lambda_{\sigma_1}}{l+1}}}^{{\Lambda_{\sigma_2}}V}(z)+
  \omega^{-\frac{1}{2}}
{{\widehat{\Psi}}_{{\Lambda_{\sigma_1}}l}}^{{\Lambda_{\sigma_2}}V}(z)
  = {{\widehat{\Psi}}_{{\Lambda_{\sigma_1}}{l+1}}}^{{\Lambda_{\sigma_2}}V}(z)
f_l$$
  $$q f_l
{{\widehat{\Psi}}_{{\Lambda_{\sigma_1}}{l+2}}}^{{\Lambda_{\sigma_2}}V}(z)+
  \omega^{-\frac{1}{2}}
{{\widehat{\Psi}}_{{\Lambda_{\sigma_1}}{l+1}}}^{{\Lambda_{\sigma_2}}V}(z)
  = {{\widehat{\Psi}}_{{\Lambda_{\sigma_1}}{l+2}}}^{{\Lambda_{\sigma_2}}V}(z)
f_l$$
  $$q f_{l-1}
{{\widehat{\Psi}}_{{\Lambda_{\sigma_1}}{l+3}}}^{{\Lambda_{\sigma_2}}V}(z)+
  {{\widehat{\Psi}}_{{\Lambda_{\sigma_1}}{l+2}}}^{{\Lambda_{\sigma_2}}V}(z)
  = {{\widehat{\Psi}}_{{\Lambda_{\sigma_1}}{l+3}}}^{{\Lambda_{\sigma_2}}V}(z)
f_{l-1}$$
  $$\vdots~~~~\vdots~~~~\vdots $$
  $$q f_2
{{\widehat{\Psi}}_{{\Lambda_{\sigma_1}}{2l}}}^{{\Lambda_{\sigma_2}}V}(z)+
  {{\widehat{\Psi}}_{{\Lambda_{\sigma_1}}{2l-1}}}^{{\Lambda_{\sigma_2}}V}(z)=
  {{\widehat{\Psi}}_{{\Lambda_{\sigma_1}}{2l}}}^{{\Lambda_{\sigma_2}}V}(z)
f_2$$
  $$q f_1
{{\widehat{\Psi}}_{{\Lambda_{\sigma_1}}{2l+1}}}^{{\Lambda_{\sigma_2}}V}(z)+
  {{\widehat{\Psi}}_{{\Lambda_{\sigma_1}}{2l}}}^{{\Lambda_{\sigma_2}}V}(z)=
  {{\widehat{\Psi}}_{{\Lambda_{\sigma_1}}{2l+1}}}^{{\Lambda_{\sigma_2}}V}(z)
f_1$$
  $$q^{-1} e^{\prime}_0
{{\widehat{\Psi}}_{{\Lambda_{\sigma_1}}{2l+1}}}^{{\Lambda_{\sigma_2}}V}(z)
  +q^2 z {{\widehat{\Psi}}_{{\Lambda_{\sigma_1}}2}}^{{\Lambda_{\sigma_2}}V}(z)=
  {{\widehat{\Psi}}_{{\Lambda_{\sigma_1}}{2l+1}}}^{{\Lambda_{\sigma_2}}V}(z)
e^{\prime}_0$$
  $$[e_i,
{{\widehat{\Psi}}_{{\Lambda_{\sigma_1}}{2l+1}}}^{{\Lambda_{\sigma_2}}V}(z)]=0,
  ~~~i=1,\cdots ,l$$

$$[e_i,{{\widehat{\Psi}}_{{\Lambda_{\sigma_1}}{i+1}}}^{{\Lambda_{\sigma_2}}V}(z)]=0,~~~i\not= l$$

$$[e_i,{{\widehat{\Psi}}_{{\Lambda_{\sigma_1}}{2l-i+1}}}^{{\Lambda_{\sigma_2}}V}(z)]=0,~~~i\not= l$$
  $$q^{h_i}
{{\widehat{\Psi}}_{{\Lambda_{\sigma_1}}{i+1}}}^{{\Lambda_{\sigma_2}}V}(z)
q^{-h_i}=q
{{\widehat{\Psi}}_{{\Lambda_{\sigma_1}}{i+1}}}^{{\Lambda_{\sigma_2}}V}(z),~~~i\not= l$$
  $$q^{h_i}
{{\widehat{\Psi}}_{{\Lambda_{\sigma_1}}{2l-i+2}}}^{{\Lambda_{\sigma_2}}V}(z)
q^{-h_i}=q
{{\widehat{\Psi}}_{{\Lambda_{\sigma_1}}{2l-i+2}}}^{{\Lambda_{\sigma_2}}V}(z),~~~i\not= l$$

$$[e_l,{{\widehat{\Psi}}_{{\Lambda_{\sigma_1}}{l+2}}}^{{\Lambda_{\sigma_2}}V}(z)]=0,
  ~~~[q_l^{h_l},
{{\widehat{\Psi}}_{{\Lambda_{\sigma_1}}{l+1}}}^{{\Lambda_{\sigma_2}}V}(z)]=0,
  ~~~q_l^{h_l}
{{\widehat{\Psi}}_{{\Lambda_{\sigma_1}}{l+2}}}^{{\Lambda_{\sigma_2}}V}(z)
q_l^{-h_l}
  =q
{{\widehat{\Psi}}_{{\Lambda_{\sigma_1}}{l+2}}}^{{\Lambda_{\sigma_2}}V}(z)$$
  $$q^{h_0}
{{\widehat{\Psi}}_{{\Lambda_{\sigma_1}}{2l+1}}}^{{\Lambda_{\sigma_2}}V}(z)
q^{-h_0}
  =q^{-1}
{{\widehat{\Psi}}_{{\Lambda_{\sigma_1}}{2l+1}}}^{{\Lambda_{\sigma_2}}V}(z)$$
Define operators ${P^{\pm}}_i$
$${P^{\pm }}_i x=[{x^{\pm}}_i(0),x] q^{\mp h_i}$$
we have $${P^{\pm }}_i {W_j}^{\mp }x=\frac{(q^{\pm h_j} x q^{\mp h_j}
-q^{\mp h_j} x q^{\pm h_j})}{q_i-{q_i}^{-1}} {\delta}_{i,j} +{W_j}^{\mp }
{P^{\pm }}_i x$$
Using them, we get
$${P^{+}}_2 {P^{+}}_3\cdots {P^{+}}_l {P^{+}}_l\cdots {P^{+}}_3
{P^{+}}_2 e^{\prime}_0
=(q^{\frac{1}{2}}+q^{-\frac{1}{2}})^2 {x^{-}}_1(1)$$
and the consistent identity
$$q^{2l-2} {\omega}
({{\widehat{\Psi}}_{{\Lambda_{\sigma_1}}{2l+1}}}^{{\Lambda_{\sigma_2}}V}(z) f_1
-q
f_1{{\widehat{\Psi}}_{{\Lambda_{\sigma_1}}{2l+1}}}^{{\Lambda_{\sigma_2}}V}(z))=
(qz)^{-1}({{\widehat{\Psi}}_{{\Lambda_{\sigma_1}}{2l+1}}}^{{\Lambda_{\sigma_2}}V}(z)
{x^{-}}_1(1)-q^{-1} {x^{-}}_1(1)
{{\widehat{\Psi}}_{{\Lambda_{\sigma_1}}{2l+1}}}^{{\Lambda_{\sigma_2}}V}(z))$$
Commuting above identity with $X_1^{+}(w)$ ,we get
$$\begin{array}{c}
\frac{q^{2l}z
{\omega}}{w}({{\widehat{\Psi}}_{{\Lambda_{\sigma_1}}{2l+1}}}^{{\Lambda_{\sigma_2}}V}(z)
{{\widehat{\varphi}}^{+}}_1(w {\gamma}^{-\frac{1}{2}}) -q^2
{{\widehat{\varphi}}^{+}}_1(w {\gamma}^{-\frac{1}{2}} )
{{\widehat{\Psi}}_{{\Lambda_{\sigma_1}}{2l+1}}}^{{\Lambda_{\sigma_2}}V}(z))\\
={{\widehat{\Psi}}_{{\Lambda_{\sigma_1}}{2l+1}}}^{{\Lambda_{\sigma_2}}V}(z)
{{\widehat{\varphi}}^{+}}_1(w {\gamma}^{-\frac{1}{2}})-
{{\widehat{\varphi}}^{+}}_1(w {\gamma}^{-\frac{1}{2}})
{{\widehat{\Psi}}_{{\Lambda_{\sigma_1}}{2l+1}}}^{{\Lambda_{\sigma_2}}V}(z)
\end{array}$$
$$\begin{array}{c}
\frac{q^{2l-2}
z{\omega}}{w}({{\widehat{\Psi}}_{{\Lambda_{\sigma_1}}{2l+1}}}^{{\Lambda_{\sigma_2}}V}(z)
{{\widehat{\varphi}}^{-}}_1(w {\gamma}^{\frac{1}{2}})-
{{\widehat{\varphi}}^{-}}_1(w {\gamma}^{\frac{1}{2}})
{{\widehat{\Psi}}_{{\Lambda_{\sigma_1}}{2l+1}}}^{{\Lambda_{\sigma_2}}V}(z))\\
={{\widehat{\Psi}}_{{\Lambda_{\sigma_1}}{2l+1}}}^{{\Lambda_{\sigma_2}}V}(z)
{{\widehat{\varphi}}^{-}}_1(w {\gamma}^{\frac{1}{2}})-q^{-2}
{{\widehat{\varphi}}^{-}}_1(w {\gamma}^{\frac{1}{2}})
{{\widehat{\Psi}}_{{\Lambda_{\sigma_1}}{2l+1}}}^{{\Lambda_{\sigma_2}}V}(z)
\end{array}$$
where ${{\widehat{\varphi}}^{\pm}}_i(w)=q^{\mp h_i} {\varphi}^{\pm}_i(w)$.
The above two formulae hold to any powers of $w$.Thus we obtain
$$[a_1(n),{{\widehat{\Psi}}_{{\Lambda_{\sigma_1}}{2l+1}}}^{{\Lambda_{\sigma_2}}V}(z)]
=\frac{[n]}{n} (q^{2l+\frac{1}{2}} z {\omega})^n
{{\widehat{\Psi}}_{{\Lambda_{\sigma_1}}{2l+1}}}^{{\Lambda_{\sigma_2}}V}(z)$$
$$[a_1(-n),{{\widehat{\Psi}}_{{\Lambda_{\sigma_1}}{2l+1}}}^{{\Lambda_{\sigma_2}}V}(z)]
=-\frac{[n]}{n} (q^{2l-\frac{1}{2}} z {\omega})^{-n}
{{\widehat{\Psi}}_{{\Lambda_{\sigma_1}}{2l+1}}}^{{\Lambda_{\sigma_2}}V}(z)$$
where $n \geq 0$.
Meanwhile from the intertwining relation, we also get
$$[e_i
,{{\widehat{\Psi}}_{{\Lambda_{\sigma_1}}{2l+1}}}^{{\Lambda_{\sigma_2}}V}(z)]=0,
{}~~~i\not= 0$$
$$[f_i
,{{\widehat{\Psi}}_{{\Lambda_{\sigma_1}}{2l+1}}}^{{\Lambda_{\sigma_2}}V}(z)]=0,~~~i\not= 1$$
Commuting the above two identities with $a_1(\pm n)$.We can easily prove
$$[x_2^{+}(w),{{\widehat{\Psi}}_{{\Lambda_{\sigma_1}}{2l+1}}}^{{\Lambda_{\sigma_2}}V}(z)]=0$$
$$[x_2^{-}(w),{{\widehat{\Psi}}_{{\Lambda_{\sigma_1}}{2l+1}}}^{{\Lambda_{\sigma_2}}V}(z)]=0$$
So we obtain
$$[a_2(n),{{\widehat{\Psi}}_{{\Lambda_{\sigma_1}}{2l+1}}}^{{\Lambda_{\sigma_2}}V}(z)]=0$$
Then repeat above procedures replacing $a_1(\pm n)$ by $a_2(\pm n)$,and so
on,we get
$$[x_i^{+}(w),{{\widehat{\Psi}}_{{\Lambda_{\sigma_1}}{2l+1}}}^{{\Lambda_{\sigma_2}}V}(z)]=0$$
$$[x_i^{-}(w),{{\widehat{\Psi}}_{{\Lambda_{\sigma_1}}{2l+1}}}^{{\Lambda_{\sigma_2}}V}(z)]=0,~~~i\not=1$$
We can conclude that
$$[\Psi (w),
{{\widehat{\Psi}}_{{\Lambda_{\sigma_1}}{2l+1}}}^{{\Lambda_{\sigma_2}}V}(z)]=0$$
$$[a_i(\pm
n),{{\widehat{\Psi}}_{{\Lambda_{\sigma_1}}{2l+1}}}^{{\Lambda_{\sigma_2}}V}(z)]=0,~~~i\not= 1$$
Put
$${a_1}^*(k)=\sum_{n=1}^{l-1} \frac{([(l-n)k]-[(l-n-1)k])}{([lk]-[(l-1)k])[k]}
a_n(k)
+\frac{[k]}{([lk]-[(l-1)k])[2k]} a_l(k)$$
we get
$$[a_i(k),{a_1}^*(-k)]={\delta}_{i,1} \frac{[k]}{k},~~~k>0$$
$$[a_i(-k),{a_1}^*(k)]=-{\delta}_{i,1} \frac{[k]}{k},~~~k<0$$
{}From the following relations:
$$[\Psi (w),
{{\widehat{\Psi}}_{{\Lambda_{\sigma_1}}{2l+1}}}^{{\Lambda_{\sigma_2}}V}(z)]=0$$
$$[a_i(\pm
n),{{\widehat{\Psi}}_{{\Lambda_{\sigma_1}}{2l+1}}}^{{\Lambda_{\sigma_2}}V}(z)]=0,~~~i\not= 1$$
$$[a_1(n),{{\widehat{\Psi}}_{{\Lambda_{\sigma_1}}{2l+1}}}^{{\Lambda_{\sigma_2}}V}(z)]
=\frac{[n]}{n} (q^{2l+\frac{1}{2}} z {\omega})^n
{{\widehat{\Psi}}_{{\Lambda_{\sigma_1}}{2l+1}}}^{{\Lambda_{\sigma_2}}V}(z),~~~n\in {\bf Z}_{>0}$$
$$[a_1(-n),{{\widehat{\Psi}}_{{\Lambda_{\sigma_1}}{2l+1}}}^{{\Lambda_{\sigma_2}}V}(z)]
=-\frac{[n]}{n} (q^{2l-\frac{1}{2}} z {\omega})^{-n}
{{\widehat{\Psi}}_{{\Lambda_{\sigma_1}}{2l+1}}}^{{\Lambda_{\sigma_2}}V}(z),~~~n\in {\bf Z}_{>0}$$
$$[x_i^{+}(w),{{\widehat{\Psi}}_{{\Lambda_{\sigma_1}}{2l+1}}}^{{\Lambda_{\sigma_2}}V}(z)]=0$$
$$q^{h_1}{{\widehat{\Psi}}_{{\Lambda_{\sigma_1}}{2l+1}}}^{{\Lambda_{\sigma_2}}V}(z) q^{-h_1}
=q {{\widehat{\Psi}}_{{\Lambda_{\sigma_1}}{2l+1}}}^{{\Lambda_{\sigma_2}}V}(z)$$
$$q^{h_i}{{\widehat{\Psi}}_{{\Lambda_{\sigma_1}}{2l+1}}}^{{\Lambda_{\sigma_2}}V}(z) q^{-h_i}
={{\widehat{\Psi}}_{{\Lambda_{\sigma_1}}{2l+1}}}^{{\Lambda_{\sigma_2}}V}(z),~~~i\geq 2 $$
$$q^{d}{{\widehat{\Psi}}_{{\Lambda_{\sigma_1}}{2l+1}}}^{{\Lambda_{\sigma_2}}V}(z) q^{-d}
={{\widehat{\Psi}}_{{\Lambda_{\sigma_1}}{2l+1}}}^{{\Lambda_{\sigma_2}}V}(q^{-1}z) $$
we obtain the type-{\rm I} vertex operators
$${{\widehat{\Psi}}_{{\Lambda_{i}}{2l+1}}}^{{\Lambda_{1-i}}V}(z)=
\omega^{-i} exp[\sum_{n=1}^\infty(q^{2l+\frac{1}{2}} z {\omega})^n {a_1}^*(-n)]
exp[-\sum_{n=1}^\infty(q^{2l-\frac{1}{2}} z {\omega})^{-n} {a_1}^* (n)]
e^{{\lambda_1}}
(q^{2l} z {\omega})^{{\partial}_{\lambda_1}+i} (-1)^{2 \partial_{\lambda_l}} $$
where $i=0,1$
$${{\widehat{\Psi}}_{{\Lambda_{l}}{2l+1}}}^{{\Lambda_{l}}V}(z)=
\omega^{-\frac{1}{2}} exp[\sum_{n=1}^\infty (q^{2l+\frac{1}{2}} z {\omega})^n
{a_1}^*(-n)]
exp[-\sum_{n=1}^\infty (q^{2l-\frac{1}{2}} z {\omega})^{-n} {a_1}^* (n)]
e^{{\lambda_1}}
(q^{2l} z {\omega})^{{\partial}_{\lambda_1}} (-1)^{2\partial_{\lambda_l}}$$
$${{\widehat{\Psi}}_{{\Lambda_{\sigma_1}}2l+1-n}}^{{\Lambda_{\sigma_2}}V}(z)
= {{\widehat{\Psi}}_{{\Lambda_{\sigma_1}}{2l+2-n}}}^{{\Lambda_{\sigma_2}}V}(z)
f_n-q f_n
{{\widehat{\Psi}}_{{\Lambda_{\sigma_1}}{2l+2-n}}}^{{\Lambda_{\sigma_2}}V}(z),
{}~~~n<l$$
$${{\widehat{\Psi}}_{{\Lambda_{\sigma_1}}n}}^{{\Lambda_{\sigma_2}}V}(z)
= {{\widehat{\Psi}}_{{\Lambda_{\sigma_1}}{n+1}}}^{{\Lambda_{\sigma_2}}V}(z)
f_n-q f_n
{{\widehat{\Psi}}_{{\Lambda_{\sigma_1}}{n+1}}}^{{\Lambda_{\sigma_2}}V}(z),
{}~~~n<l$$
$${{\widehat{\Psi}}_{{\Lambda_{\sigma_1}}{l+1}}}^{{\Lambda_{\sigma_2}}V}(z)=
{\omega}^{\frac{1}{2}}
({{\widehat{\Psi}}_{{\Lambda_{\sigma_1}}{l+2}}}^{{\Lambda_{\sigma_2}}V}(z) f_l-
q f_l
{{\widehat{\Psi}}_{{\Lambda_{\sigma_1}}{l+2}}}^{{\Lambda_{\sigma_2}}V}(z))$$
$${{\widehat{\Psi}}_{{\Lambda_{\sigma_1}}l}}^{{\Lambda_{\sigma_2}}V}(z)
={\omega}^{\frac{1}{2}}(
{{\widehat{\Psi}}_{{\Lambda_{\sigma_1}}{l+1}}}^{{\Lambda_{\sigma_2}}V}(z) f_l
-f_l {{\widehat{\Psi}}_{{\Lambda_{\sigma_1}}{l+1}}}^{{\Lambda_{\sigma_2}}V}(z)
)$$
where we have used the nomalization condition for $i=0,1,l$.

\noindent
$3.4$ Vertex Operator of Type-{\rm II}

For the Vertex Operator of type-{\rm II},the analogous commutation relations
can be got.
We have
$$[{\Psi}(W),{{\widehat{\Phi}}_{{\Lambda_{\sigma_2}}1}}^{V{\Lambda_{\sigma_1}}}(z)]=0$$
$$[a_1(n),{{\widehat{\Phi}}_{{\Lambda_{\sigma_2}}1}}^{V{\Lambda_{\sigma_1}}}(z)]
=-\frac{[n]}{n} (q^{\frac{1}{2}} {\omega} z)^n
{{\widehat{\Phi}}_{{\Lambda_{\sigma_2}}1}}^{V{\Lambda_{\sigma_1}}}(z)$$
$$[a_1(-n),{{\widehat{\Phi}}_{{\Lambda_{\sigma_2}}1}}^{V{\Lambda_{\sigma_1}}}(z)]
=\frac{[n]}{n} (q^{-\frac{1}{2}} {\omega} z)^{-n}
{{\widehat{\Phi}}_{{\Lambda_{\sigma_2}}1}}^{V{\Lambda_{\sigma_1}}}(z)$$
$$[a_i(\pm
n),{{\widehat{\Phi}}_{{\Lambda_{\sigma_2}}1}}^{V{\Lambda_{\sigma_1}}}(z)]=0,
{}~~~i\not= 1$$
$$[x_i^{-},{{\widehat{\Phi}}_{{\Lambda_{\sigma_2}}1}}^{V{\Lambda_{\sigma_1}}}(z)]=0$$
$$q^{h_1} {{\widehat{\Phi}}_{{\Lambda_{\sigma_2}}1}}^{V{\Lambda_{\sigma_1}}}(z)
q^{-h_1}
=q^{-1} {{\widehat{\Phi}}_{{\Lambda_i}1}}^{V{\Lambda_{\sigma_1}}}(z)$$
$$q^{h_i} {{\widehat{\Phi}}_{{\Lambda_{\sigma_2}}1}}^{V{\Lambda_{\sigma_1}}}(z)
q^{-h_i}=
{{\widehat{\Phi}}_{{\Lambda_{\sigma_2}}1}}^{V{\Lambda_{\sigma_1}}}(z),~~~i\geq
2$$
$$q^{-d} {{\widehat{\Phi}}_{{\Lambda_{\sigma_2}}1}}^{V{\Lambda_{\sigma_1}}}(z)
q^d=
{{\widehat{\Phi}}_{{\Lambda_{\sigma_2}}1}}^{V{\Lambda_{\sigma_1}}}(qz)$$
So we can get the type-{\rm II} vertex operators
$${{\widehat{\Phi}}_{{\Lambda_i}1}}^{V{\Lambda_{1-i}}}(z)=
( \omega q^{2l-1})^{i-1} exp[-\sum_{n=1}^\infty(q^{\frac{1}{2}} z {\omega})^n
{a_1}^*(-n)]
exp[\sum_{n=1}^\infty(q^{-\frac{1}{2}} z {\omega})^{-n} {a_1}^* (n)]
e^{-{\lambda_1}}
(z {\omega})^{-{\partial}_{\lambda_1}+i}(-1)^{2 \partial_{\lambda_l}}$$
where $i=0,1$
$${{\widehat{\Phi}}_{{\Lambda_l}1}}^{V{\Lambda_l}}(z)=
\omega^{-\frac{1}{2}} (-q)^{-l} exp[-\sum_{n=1}^\infty(q^{\frac{1}{2}} z
{\omega})^n {a_1}^*(-n)]
exp[\sum_{n=1}^\infty(q^{-\frac{1}{2}} z {\omega})^{-n} {a_1}^* (n)]
e^{-{\lambda_1}}
(z {\omega})^{-{\partial}_{\lambda_1}}(-1)^{2 \partial_{\lambda_l}}$$
$${{\widehat{\Phi}}_{{\Lambda_{\sigma_2}}{n+1}}}^{V{\Lambda_{\sigma_1}}}(z)=
{{\widehat{\Phi}}_{{\Lambda_{\sigma_2}}{n}}}^{V{\Lambda_{\sigma_1}}}(z) e_n-
q e_n {{\widehat{\Phi}}_{{\Lambda_i}{n}}}^{V{\Lambda_j}}(z),
{}~~~n<l$$
$${{\widehat{\Phi}}_{{\Lambda_{\sigma_2}}{2l+2-n}}}^{V{\Lambda_{\sigma_1}}}(z)=
{{\widehat{\Phi}}_{{\Lambda_{\sigma_2}}{2l+1-n}}}^{V{\Lambda_{\sigma_1}}}(z)
e_n-
q e_n {{\widehat{\Phi}}_{{\Lambda_i}{2l+1-n}}}^{V{\Lambda_j}}(z),
{}~~~n<l$$
$${{\widehat{\Phi}}_{{\Lambda_{\sigma_2}}{l+2}}}^{V{\Lambda_{\sigma_1}}}(z)=
{\omega}^{\frac{1}{2}}
({{\widehat{\Phi}}_{{\Lambda_{\sigma_2}}{l+1}}}^{V{\Lambda_{\sigma_1}}}(z) e_l-
e_l
{{\widehat{\Phi}}_{{\Lambda_{\sigma_2}}{{l+1}}}}^{V{\Lambda_{\sigma_1}}}(z))$$
$${{\widehat{\Phi}}_{{\Lambda_{\sigma_2}}{l+1}}}^{V{\Lambda_{\sigma_1}}}(z)={\omega}^{\frac{1}{2}}
( {{\widehat{\Phi}}_{{\Lambda_{\sigma_2}}l}}^{V{\Lambda_{\sigma_1}}}(z) e_l-q
e_l {{\widehat{\Phi}}_{{\Lambda_{\sigma_2}}l}}^{V{\Lambda_{\sigma_1}}}(z) )$$

\noindent
$3.5$ Dual Vertex Operators

\noindent
Define the intertwiners of the form
$${{\widehat{\Psi^*}}^{\Lambda_i}}_{{\Lambda_j}V}(z) :
V(\Lambda_j) \otimes V_z\rightarrow  V(\Lambda_i)$$
$${{\widehat{\Phi^*}}^{\Lambda_i}}_{V{\Lambda_j}}(z) :
V_z\otimes V(\Lambda_j) \rightarrow  V(\Lambda_i)$$
They are called dual vertex operators.Define their components by
$${{\widehat{\Psi^*}}^{\Lambda_i}}_{{\Lambda_j}Vm}(z) \mid v \rangle
{}~=~{{\widehat{\Psi^*}}^{\Lambda_i}}_{{\Lambda_j}V}(z) (\mid v \rangle
\otimes v_m )$$
$${{\widehat{\Phi^*}}^{\Lambda_i}}_{V{\Lambda_j}m}(z) \mid v \rangle
{}~=~{{\widehat{\Psi^*}}^{\Lambda_i}}_{V{\Lambda_j}}(z) (v_m \otimes
\mid v \rangle )$$
We impose the normalization
$$\langle \Lambda_1 \mid {{\widehat{\Psi^*}}^{\Lambda_1}}_{{\Lambda_0}V1}
(z)\mid \Lambda_0 \rangle ~=~ 1,~~~
\langle \Lambda_1 \mid  {{\widehat{\Phi^*}}^{\Lambda_1}}_{V{\Lambda_0}1}
(z) \mid \Lambda_0 \rangle ~=~ 1$$
$$\langle \Lambda_0 \mid
{{\widehat{\Psi^*}}^{\Lambda_0}}_{{\Lambda_1}V{2l+1}}(z)
\mid \Lambda_1 \rangle ~=~ 1,~~~
\langle \Lambda_0 \mid {{\widehat{\Phi^*}}^{\Lambda_0}}_{V{\Lambda_1}{2l+1}}(z)
\mid \Lambda_1 \rangle ~=~ 1$$
$$\langle \Lambda_l^{\prime} \mid
{{\widehat{\Psi^*}}^{\Lambda_l^{\prime}}}_{{\Lambda_l}V{l+1}}(z) \mid
\Lambda_l \rangle~=~1,~~~
\langle \Lambda_l^{\prime} \mid
{{\widehat{\Phi^*}}^{\Lambda_l^{\prime}}}_{V{\Lambda_l}{l+1}}(z)\mid
\Lambda_l \rangle ~=~ 1$$
With analogue dicussion in \cite{JM} p79, we get
$${{\widehat{\Phi^*}}^{\Lambda_{1-i}}}_{V{\Lambda_i}{2l+1}}(z)
{}~=~q^{(2l-1)(1-i)} {{\widehat{\Phi}}_{\Lambda_i
1}}^{V{\Lambda_{1-i}}}(q^{2l-1} z)$$
$${{\widehat{\Psi^*}}^{\Lambda_{1-i}}}_{{\Lambda_i}V1}(z)
{}~=~q^{(2l-1) i} {{\widehat{\Psi}}_{\Lambda_i
{2l+1}}}^{{\Lambda_{1-i}}V}(q^{-2l+1} z) $$
$${{\widehat{\Phi^*}}^{\Lambda_l}}_{V{\Lambda_l}{2l+1}}(z)
{}~=~(-q)^l {{\widehat{\Phi}}_{\Lambda_l 1}}^{V{\Lambda_l}}(q^{2l-1} z)$$
$${{\widehat{\Psi^*}}^{\Lambda_l}}_{{\Lambda_l}V1}(z)
{}~=~(-q)^l {{\widehat{\Psi}}_{\Lambda_l{2l+1}}}^{{\Lambda_l}V}(q^{-2l+1} z) $$
$${{\widehat{\Psi^*}}_{{\Lambda_{\sigma_1}}2l+2-n}}^{{\Lambda_{\sigma_2}}V}(z)
=f_n
{{\widehat{\Psi^*}}_{{\Lambda_{\sigma_1}}{2l+1-n}}}^{{\Lambda_{\sigma_2}}V}(z)
-q^{-1}
{{\widehat{\Psi^*}}_{{\Lambda_{\sigma_1}}{2l+1-n}}}^{{\Lambda_{\sigma_2}}V}(z)
f_n, ~~~n<l$$
$${{\widehat{\Psi^*}}_{{\Lambda_{\sigma_1}}{n+1}}}^{{\Lambda_{\sigma_2}}V}(z)
=f_n {{\widehat{\Psi^*}}_{{\Lambda_{\sigma_1}}n}}^{{\Lambda_{\sigma_2}}V}(z)
-q^{-1}
{{\widehat{\Psi^*}}_{{\Lambda_{\sigma_1}}n}}^{{\Lambda_{\sigma_2}}V}(z) f_n,
{}~~~n<l$$
$${{\widehat{\Psi^*}}_{{\Lambda_{\sigma_1}}{l+2}}}^{{\Lambda_{\sigma_2}}V}(z)=
{\omega^{\frac{1}{2}}} (f_l
{{\widehat{\Psi^*}}_{{\Lambda_{\sigma_1}}{l+1}}}^{{\Lambda_{\sigma_2}}V}(z) -
{{\widehat{\Psi^*}}_{{\Lambda_{\sigma_1}}{l+1}}}^{{\Lambda_{\sigma_2}}V}(z)
f_l)$$
$${{\widehat{\Psi^*}}_{{\Lambda_{\sigma_1}}{l+1}}}^{{\Lambda_{\sigma_2}}V}(z)
={\omega^{\frac{1}{2}}}(f_l
{{\widehat{\Psi^*}}_{{\Lambda_{\sigma_1}}l}}^{{\Lambda_{\sigma_2}}V}(z)
-q^{-1} {{\widehat{\Psi^*}}_{{\Lambda_{\sigma_1}}l}}^{{\Lambda_{\sigma_2}}V}(z)
f_l )$$
$${{\widehat{\Phi^*}}_{{\Lambda_{\sigma_2}}n}}^{V{\Lambda_{\sigma_1}}}(z)=
e_n {{\widehat{\Phi^*}}_{{\Lambda_{\sigma_2}}{n+1}}}^{V{\Lambda_{\sigma_1}}}(z)
-
q^{-1} {{\widehat{\Phi^*}}_{{\Lambda_i}{n+1}}}^{V{\Lambda_j}}(z) e_n,
{}~~~n<l$$
$${{\widehat{\Phi^*}}_{{\Lambda_{\sigma_2}}{2l+1-n}}}^{V{\Lambda_{\sigma_1}}}(z)=
e_n
{{\widehat{\Phi^*}}_{{\Lambda_{\sigma_2}}{2l+2-n}}}^{V{\Lambda_{\sigma_1}}}(z)
-
q^{-1} {{\widehat{\Phi^*}}_{{\Lambda_i}{2l+2-n}}}^{V{\Lambda_j}}(z) e_n,
{}~~~n<l$$
$${{\widehat{\Phi^*}}_{{\Lambda_{\sigma_2}}{l+1}}}^{V{\Lambda_{\sigma_1}}}(z)=
 {\omega^{\frac{1}{2}}} (e_l
{{\widehat{\Phi^*}}_{{\Lambda_{\sigma_2}}{l+2}}}^{V{\Lambda_{\sigma_1}}}(z) -
q^{-1}
{{\widehat{\Phi^*}}_{{\Lambda_{\sigma_2}}{l+2}}}^{V{\Lambda_{\sigma_1}}}(z)
e_l)$$
$${{\widehat{\Phi^*}}_{{\Lambda_{\sigma_2}}l}}^{V{\Lambda_{\sigma_1}}}(z)={\omega^{\frac{1}{2}}}
(e_l
{{\widehat{\Phi^*}}_{{\Lambda_{\sigma_2}}{l+1}}}^{V{\Lambda_{\sigma_1}}}(z) -
 {{\widehat{\Phi^*}}_{{\Lambda_{\sigma_2}}{l+1}}}^{V{\Lambda_{\sigma_1}}}(z)e_l
)$$

\noindent
{\bf Acknowledgments}. This work was partly supported by the Foundation
of the University of Science and Tecnology of China and the National
Natural Science Foundation of China and Grant No. LWTZ-1298 of
Chinese Academy of Sciences. The author[B.Jin] also thanks Prof. Hou Bo-Yu
for focusing his attention on this problem.

\end{document}